\documentclass[12pt,oneside,reqno]{amsart}
\usepackage{graphicx}
\usepackage{mathrsfs}
\usepackage{amsfonts}
\usepackage{enumerate,amsmath,amssymb,amsthm}
\newcommand{\dif}{\mathrm{d}}
\newcommand{\me}{\mathrm{e}}
\newcommand{\be}{\begin{eqnarray}}
\newcommand{\ee}{\end{eqnarray}}
\newcommand{\ce}{\begin{eqnarray*}}
\newcommand{\de}{\end{eqnarray*}}
\newtheorem{theorem}{Theorem}[section]
\newtheorem{lemma}[theorem]{Lemma}
\newtheorem{remark}[theorem]{Remark}
\newtheorem{definition}[theorem]{Definition}
\newtheorem{proposition}[theorem]{Proposition}
\newtheorem{Examples}[theorem]{Examples}
\newtheorem{corollary}[theorem]{Corollary}

\def\a{\alpha}

\def\[{{\Big[}}
\def\]{{\Big]}}
\def\<{{\langle}}
\def\>{{\rangle}}
\def\({{\Big(}}
\def\){{\Big)}}

\def\no{\nonumber}
\def\bt{\begin{theorem}}
\def\et{\end{theorem}}
\def\bl{\begin{lemma}}
\def\el{\end{lemma}}
\def\br{\begin{remark}}
\def\er{\end{remark}}
\def\bx{\begin{Examples}}
\def\ex{\end{Examples}}
\def\bd{\begin{definition}}
\def\ed{\end{definition}}
\def\bp{\begin{proposition}}
\def\ep{\end{proposition}}
\def\bc{\begin{corollary}}
\def\ec{\end{corollary}}

\def\cB{{\mathcal B}}
\def\cC{{\mathcal C}}
\def\cD{{\mathcal D}}

\def\cH{{\mathcal H}}

\def\cL{{\mathcal L}}
\def\cM{{\mathcal M}}

\def\cU{{\mathcal U}}

\def\mE{{\mathbb E}}

\def\mH{{\mathbb H}}

\def\mP{{\mathbb P}}

\def\mR{{\mathbb R}}

\def\mU{{\mathbb U}}

\def\mX{{\mathbb X}}

\def\sB{{\mathscr B}}

\def\sF{{\mathscr F}}

\def\geq{\geqslant}
\def\leq{\leqslant}

\begin{document}

\allowdisplaybreaks

\title{Stationary Measures for Stochastic Differential Equations with Jumps*}

\author{Huijie Qiao and Jinqiao Duan}



\thanks{*This work was partly supported by the NSF of China (No. 11001051).}


\date{}

\dedicatory{Department of Mathematics,
Southeast University,\\
Nanjing, Jiangsu 211189, China\\
hjqiaogean@seu.edu.cn\\
Department of Applied Mathematics, Illinois Institute of Technology,\\
Chicago, IL 60616, USA\\
duan@iit.edu}

\begin{abstract}
In the paper, stationary measures of stochastic differential equations with jumps are considered.
Under some general conditions, existence of stationary measures is proved through Markov measures and
Lyapunov functions. Moreover, for two special cases, stationary measures are given by solutions
of Fokker-Planck equations and long time limits for the distributions of system states.
\end{abstract}

\maketitle \rm

\section{Introduction}
Stationary measures for stochastic differential equations (SDEs) and invariant measures for Markov processes  have been studied
extensively,  as in \cite{da, la, hc, cs, jz}, among others. Moreover, there are interesting relationships between stationary
measures and invariant measures,  such as the one-to-one correspondence between the set of invariant Markov measures
(Definition \ref{Marm} in Section \ref{exst}) and the set of stationary measures (\cite{la, hc}), as well as the correspondence
between the solutions of Fokker-Planck equations and stationary measures (\cite{cs}). These results play an important role in
the development of the theory for random dynamical systems associated with SDEs.

Existence and uniqueness of the stationary measures of a one dimensional diffusion process with Gaussian noise
have been discussed in \cite{wata}. Khasminski \cite{khas} extended the results in not only for higher dimension
but also for a Markov process framework so that any solution process of a parabolic SDE (with certain properties, e.g
strongly Feller) can be shown to have a unique stationary measure. Zakai \cite{zakai} was able to release the strong
Feller condition from Khasminski's work in Gaussian noise case and showed the existence of the stationary
measure using Lyapunov functions. Liu-Mandrekar offered a weaker condition of ultimately boundedness \cite{lima}
for construction of a Lyapunov function and thereby carried on Zakai's work further for the stationary measure
for the solution of Gaussian SDE.

In recent years stationary measures of SDEs with jumps are considered by a number of authors (\cite{arw, da, bhan, bhka, sa, jz}).
Let us mention some works. Zabczyk \cite{jz} studied stationary measures for {\it linear} SDEs with jumps.
Later, Albeverio-R\"udiger-Wu \cite{arw} discussed  stationary measures
for SDEs with jumps, in the context of  L\'evy type operators and considered mainly infinitesimal invariant measures.
The concept of infinitesimal invariant measures there is weaker than that of the usual
stationary measures. Recently, Bhan-Chankraborty-Mandrekar \cite{bhan} studied the stationary
measure for the solution of a stochastic differential equation driven by jump L\'evy processes
by the same Lyapunov approach as in \cite{miya}. In \cite{bhan} the stochastic differential equation
is $1$-dimensional and driven only by a jump L\'evy process and no Brownian motion.

We ask, naturally, whether the above correspondences for usual SDEs also hold for SDEs
with jumps. This question will be answered in Sections \ref{exst} and \ref{fpe} in the present paper. To our
knowledge, a result presented in \cite{bhka} is the closest to the result in Section \ref{fpe} of this paper.
In \cite[Theorem 3.2]{bhka} Bhatt-Karandikar set about a martingale problem and showed existence
of stationary measures for Markov processes. There, they required that the domain of the generator
is an algebra that separates points and vanishes nowhere, which is not the case in our setting.
Because the generators of SDEs with jumps which we consider here are integro-differential operators,
their domains are usually not algebras with some properties.

Now we briefly sketch our method. We begin from SDEs with jumps  as random dynamical systems
and examine their Markov measures. Then we weaken some conditions and still obtain existence
of stationary measures by the similar Lyapunov approach to in \cite{bhan}. Here, the type of
our equations is more general and these conditions are easier to satisfy than the type and
conditions in \cite{bhan}. Next, in the first special case, we consider SDEs which are driven by Brownian motions
and $\alpha$-stable processes. By a functional analysis technique, stationary measures are investigated.
When the coefficients of the SDEs are sufficiently regular in the second special case, the long time limits for the
distributions of the solutions  are shown to be stationary measures.

This paper is arranged as follows. In Section \ref{pre}, we introduce random dynamical systems and
related concepts. Symmetric $\alpha$-stable processes are also introduced. The content
to obtain existence of stationary measures through Markov measures and Lyapunov functions is in Section \ref{exst}. In
Section \ref{fpe}, we consider special stochastic differential equations
driven by Brownian motions and $\alpha$-stable processes. In Section \ref{des}, we study SDEs with jumps, under
certain regular conditions on the coefficients.

\section{Preliminary}\label{pre}

In this section, we recall  basic concepts and facts that will be
needed throughout the paper.

\subsection{Random dynamical systems and related concepts}

\bd\label{mds}
Let $(\Omega,\sF,\mP)$ be a probability space, and $(\theta_t)_{t\in\mR_+}$ a family of measurable transformations from $\Omega$ to $\Omega$.
We call $(\Omega,\sF,\mP; (\theta_t)_{t\in\mR_+})$ a metric dynamical system if for each $t\in\mR_+$, $\theta_t$ preserves the probability measure $\mP$, i.e.,
$$
\theta_t^*\mP=\mP.
$$
\ed

Throughout the paper $(\theta_t)_{t\in\mR_+}$ will be assumed ergodic, i.e. all measurable $\theta$-invariant sets have
probability either $0$ or $1$ (\cite{la}).

\bd\label{rds}
A measurable c\`adl\`ag random dynamical system on the Polish space $(\mX,\sB)$
over a metric DS $(\Omega,\mathscr{F},\mP,(\theta_t)_{t\in\mR_+})$ with time $\mR_+$ is
a family of homeomorphisms of $\mX$,
\ce
&&\varphi: \mR_+\times\Omega\times\mX\mapsto\mX, \quad
(t,\omega,x)\mapsto\varphi(t,\omega,x),\\
&&\varphi(t,\omega)\cdot:=\varphi(t,\omega,\cdot): \mX\mapsto\mX,
\de
such that

(i) Measurability: $\varphi$ is $\sB(\mR_+)\otimes\mathscr{F}\otimes\sB/\sB$-measurable,
where $\sB(\mR_+)$ is Borel $\sigma$-algebra of $\mR_+$.

(ii) C\`adl\`ag cocycle property: $\varphi(t,\omega)$ forms a (perfect) c\`adl\`ag
cocycle over $\theta$ if it is c\`adl\`ag in $t$ and satisfies
\be
\varphi(0,\omega)&=&id_{\mX}, ~\mbox{for~ all}~\omega\in\Omega,\label{perfect coc1}\\
\varphi(t+s,\omega)&=&\varphi(t,\theta_s\omega)\circ\varphi(s,\omega),
\label{perfect coc2}
\ee
for all $s,t\in\mR_+$ and $\omega\in\Omega$.
\ed

A random dynamical system (RDS) induces a skew product flow of
measurable maps
\ce
\Theta_t:\Omega\times\mX&\mapsto&\Omega\times\mX\\
(\omega,x)&\mapsto&\big(\theta_t\omega,\varphi(t,\omega)x\big).
\de
The flow property $\Theta_{t+s}=\Theta_t\circ\Theta_s$ follows from (\ref{perfect coc2}).
Denote by $\mathscr{P}(\Omega\times\mX)$ the probability measures
on $(\Omega\times\mX,\mathscr{F}\otimes\sB)$.  Moreover, $\Theta_t$ acts on $\mu\in\mathscr{P}(\Omega\times\mX)$
by $(\Theta_t\mu)(C)=\mu(\Theta_t^{-1}C)$,  for $C\in\mathscr{F}\otimes\sB$,
$t\in\mR_+$.

\bd\label{invm} A probability measure
$\mu\in\mathscr{P}(\Omega\times\mX)$ is called invariant for the
skew product flow $\Theta_t$ if

(i) the marginal of $\mu$ on $\Omega$ is $\mP$,

(ii) $\Theta_t\mu=\mu$ for all $t\in\mR_+$.
\ed

If $\mX$ is a Polish space with its Borel $\sigma$-algebra $\sB(\mX)$, every
measure $\mu\in\mathscr{P}(\Omega\times\mX)$ with marginal $\mP$ can
be uniquely characterized through its factorization
$$
\mu(\dif \omega,\dif x)=\mu_\omega(\dif x)\mP(\dif \omega),
$$
where $\mu_\omega(\dif x)$ is a probability kernel, i.e. for any $B\in\sB(\mX)$,
$\mu_{\cdot}(B)$ is $\mathscr{F}$-measurable; for $\mP.a.s.\omega\in\Omega$,
$\mu_\omega(\cdot)$ is a probability measure on $(\mX,\sB(\mX))$ (\cite[p.23]{la}).
Thus $\mu$ is invariant if and only if
\be
\mE[\varphi(t,\omega)\mu_\cdot|\theta_t^{-1}\mathscr{F}](\omega)=\mu_{\theta_t\omega},
\quad \mP. a.s,
\label{inveq}
\ee
for all $t\in\mR_+$.

\subsection{Symmetric $\alpha$-stable processes}

\bd\label{levy}
A process $L=(L_t)_{t\geq0}$ with $L_0=0$ a.s. is a
$d$-dimensional L\'evy process if

(i) $L$ has independent increments; that is, $L_t-L_s$ is
independent of $L_v-L_u$ if $(u,v)\cap(s,t)=\emptyset$;

(ii) $L$ has stationary increments; that is, $L_t-L_s$ has the same
distribution as $L_v-L_u$ if $t-s=v-u>0$;

(iii) $L_t$ is stochastically continuous;

(iv) $L_t$ is right continuous with left limit.
\ed
Its characteristic function is given by
\ce
\mE\left(\exp\{i\<z,L_t\>\}\right)=\exp\{t\Psi(z)\}, \quad
z\in\mR^d.
\de
The function $\Psi: \mR^d\rightarrow\mathcal {C}$ is
called the characteristic exponent of the L\'evy process $L$. By the
L\'evy-Khintchine formula, there exist a nonnegative-definite
$d\times d$ matrix $Q$, a measure $\nu$ on $\mR^d$ satisfying \ce
\nu(\{0\})=0 ~\mbox{and}~ \int_{\mR^d\setminus\{0\}}(|u|^2\wedge1)\nu(\dif
u)<\infty, \de and $\gamma\in\mR^d$ such that
\be
\Psi(z)=i\<z,\gamma\>-\frac{1}{2}\<z,Qz\>+\int_{\mR^d\setminus\{0\}}\big(e^{i\<z,u\>}-1-i\<z,u\>1_{|u|\leq1}\big)\nu(\dif
u).
\label{lkf}
\ee
The measure $\nu$ is called the L\'evy measure.

\bd\label{rid}
For $\alpha\in(0,2)$. A $d$-dimensional symmetric $\alpha$-stable
process $L$ is a L\'evy process such that its characteristic exponent $\Psi$
is given by for $z\in\mR^d$
\ce
\Psi(z)=-C|z|^\alpha, \quad
C=\pi^{-1/2}\frac{\Gamma((1+\a)/2)\Gamma(d/2)}{\Gamma((d+\a)/2)}.
\de
\ed
Thus, for $d$-dimensional symmetric $\alpha$-stable process $L$,
the diffusion matrix $Q=0$, the drift vector $\gamma=0$,  and the
L\'evy measure $\nu$ is given by
$$
\nu(\dif u)=\frac{C_{d,\alpha}}{|u|^{d+\alpha}}\dif u, \quad
C_{d,\alpha}=\frac{\a\Gamma((d+\a)/2)}{2^{1-\a}\pi^{d/2}\Gamma(1-\a/2)}.
$$

To give the infinitesimal generator of a $d$-dimensional symmetric $\alpha$-stable process,
we introduce several function spaces. Let $\cC_0(\mR^{d})$ be the space of continuous functions $f$ on $\mR^{d}$ satisfying
$\lim\limits_{|x|\rightarrow\infty}f(x)=0$. Let $\cC^2_0(\mR^{d})$ be the set of $f\in\cC_0(\mR^{d})$ such that $f$ is two
times differentiable and the partial derivatives of $f$ with order two and less than two belong
to $\cC_0(\mR^{d})$. Let $\cC_c^n(\mR^{d})$ stand for the space of all $n$
times differentiable functions on $\mR^{d}$ with compact supports. Let $S^\prime(\mR^{d})$
be the space of all tempered distributions on $\mR^{d}$ and $\hat{f}$ the
Fourier transform of $f\in S^\prime(\mR^{d})$. Set
\ce
\mH^{\lambda,2}(\mR^{d}):=\{f\in S^\prime(\mR^{d}):\|f\|_{\lambda,2}<\infty\},
\de
for any $\lambda\in\mR$, where
\ce
\|f\|^2_{\lambda,2}:=\int_{\mR^{d}}(1+|u|^2)^\lambda|\hat{f}(u)|^2\dif u.
\de
In particular, $\mH^{0,2}(\mR^{d})=L^2(\mR^{d})$.

Define
\ce
(\cL_\alpha h)(x):=\int_{\mR^d\setminus\{0\}}\big(h(x+u)-h(x)
-\<\partial_x h(x),u\>1_{|u|\leq1}\big)\frac{C_{d,\alpha}}{|u|^{d+\alpha}}\dif u
\de
for $h\in\cC_0^2(\mR^{d})$ and then $\cL_\alpha$ is the infinitesimal generator of a $d$-dimensional
symmetric $\alpha$-stable process(\cite{da}). Moreover, by \cite[Example 3.3.8, P.166]{da} for every $f\in\cC_c^\infty(\mR^{d})$
\ce
(\cL_\alpha f)(x)=C[-(-\Delta)^{\alpha/2}f](x).
\de
The following result from \cite{arw} is used in Section \ref{fpe}.

\bt\label{sead}
Let $\cL_\alpha$ be as above for $\alpha\in(0,2)$ and $\cL_2=\Delta$,
as defined on $\cC_c^\infty(\mR^{d})$ in $L^2(\mR^{d})$. Then $\cL_\alpha$,
$0<\alpha\leq 2$, has a unique closed extensions to self-adjoint negative
operators on the domain $\mH^{\alpha,2}(\mR^{d})$.
\et

\section{Existence of stationary measures}\label{exst}

In the section, we prove existence of stationary measures for general SDEs with jumps
under some general conditions.

Let $(\mU,\cU,n)$ be a $\sigma$-finite measurable
space. Let $\{W(t)\}_{t\geq 0}$ be an $m$-dimensional standard Brownian motion,
and $\{k_t,t\geq 0\}$ a stationary Poisson point process with values in $\mU$ and with
characteristic measure $n$ (\cite{da}). Let $N_{k}((0,t],\dif
u)$ be the counting measure of $k_{t}$, i.e., for $A\in\cU$
$$
N_k((0,t],A):=\#\{0<s\leq t: k_s\in A\},
$$
where $\#$ denotes the cardinality of a set. The compensator measure
of $N_k$ is given by
$$
\tilde{N}_{k}((0,t],\dif u):=N_{k}((0,t],\dif u)-tn(\dif u).
$$

Fix a $\mU_{0}\in\cU$ such that $n(\mU-\mU_{0})<\infty$,
and consider the following SDE with jumps in $\mR^{d}$:
\be
X_{t}(x)&=&x+\int^{t}_{0}b(X_{s}(x))\,\dif s+\int^{t}_{0}\sigma(X_{s}(x))\,\dif W_{s}\no\\
&&+\int^{t+}_{0}\int_{\mU_{0}}f(X_{s-}(x),u)\,\tilde{N}_{k}(\dif s,\dif u)\no\\
&&+\int^{t+}_{0}\int_{\mU-\mU_{0}}g(X_{s-}(x),u)\,N_{k}(\dif s,\dif u), \quad t\geq0,
\label{Eqj1}
\ee
where $b:\mR^d\mapsto\mR^d$, $\sigma:\mR^d\mapsto\mR^d\times\mR^m$,
$f,g:\mR^d\times\mU\mapsto\mR^d$ satisfy the following assumptions:
\begin{enumerate}[{\bf (H$_b$)}]
\item there exists a constant $C_b>0$ such that
for $x,y\in\mR^{d}$
$$
|b(x)-b(y)|\leq C_b|x-y|\cdot \log(|x-y|^{-1}+e);
$$
\end{enumerate}
\begin{enumerate}[{\bf (H$_\sigma$)}]
\item there exists a constant $C_\sigma>0$ such that
for $x,y\in\mR^{d}$
$$
|\sigma(x)-\sigma(y)|^2\leq C_\sigma|x-y|^{2}\cdot
\log(|x-y|^{-1}+e);
$$
\end{enumerate}

\begin{enumerate}[{\bf (H$_f$)}]
\item There exists a positive function $L(u)$ satisfying
\ce
\sup_{u\in\mU_0}L(u)\leq \delta<1~\mbox{ and } \int_{\mU_0}L(u)^2\,\nu(\dif u)<+\infty,
\de
such that for any $x,y\in\mR^{d}$ and $u\in\mU_0$
$$
|f(x,u)-f(y,u)|\leq L(u)\cdot|x-y|,
$$
and
$$
|f(0,u)|\leq L(u).
$$
Moreover, we also require that for some $q>4d$
\ce
\frac{q\delta}{(1-\delta)^{q+1}}<1.
\de
\end{enumerate}
\begin{enumerate}[{\bf (H$_g$)}]
\item For each $u\in\mU-\mU_0$, $x\mapsto x+g(x,u)\in\cH(\mR^d)$, where
$\cH(\mR^d)$ denotes the set of  all homeomorphism mappings from $\mR^d$ to $\mR^d$.
\end{enumerate}
Here, the second integral of the right hand side in Eq.(\ref{Eqj1}) is taken in It\^o's sense,
and the definitions of the third and fourth integrals are referred to \cite{da}. Under {\bf (H$_b$)},
{\bf (H$_\sigma$)}, {\bf (H$_f$)} and {\bf (H$_g$)}, it
is well known that there exists a unique strong solution to Eq.(\ref{Eqj1}) (\cite{hqxz}).
This solution will be denoted by $X_t(x)$.

Set
\ce
&&\mathscr{F}_{\geq t}:=\sigma\{W_s,N_p((0,s],B); s\geq t, B\in\cU\},\\
&&\mathscr{F}_{\leq t}:=\sigma\{W_s,N_p((0,s],B); s\leq t, B\in\cU\},
\de
for $t\geq0$. And then
$$
\theta_t^{-1}\mathscr{F}_{>0}\subset\mathscr{F}_{>t}, \quad
\theta_t^{-1}\mathscr{F}_{\leq s}\subset\mathscr{F}_{\leq t+s}, \quad s\geq0.
$$
Moreover, $X_t(x)$ is $\mathscr{F}_{\leq t}$-measurable
and independent of $\theta_t^{-1}\mathscr{F}_{>0}$.

\bd\label{Marm}
A probability measure $\mu$ on $(\Omega\times\mR^d,\mathscr{F}\otimes\sB(\mR^d))$
is called a Markov measure if $\mu_{\omega}$ satisfies
\ce
\mE(\mu_{\cdot}|\mathscr{F}_{<\infty})=\mE(\mu_{\cdot}|\mathscr{F}_{=0}), \quad \mP. a.s..
\de
\ed

For the Markov process $X_t(x)$, the transition probability is defined by
$$
p_t(x,B):=\mP\big(X_t(x)\in B\big), \quad t>0, B\in\sB(\mR^d).
$$

\bd\label{stam}
A measure $\bar{\mu}$ on $(\mR^d,\sB(\mR^d))$ is a stationary measure for $p_t$
or Eq.(\ref{Eqj1}) if
\ce
\int_{\mR^d}p_t(x,B)\bar{\mu}(\dif x)=\bar{\mu}(B), \qquad \forall t>0, B\in\sB(\mR^d).
\de
\ed

Before stating the main result in the section, we prove an important lemma.

\bl\label{invmeqstam}
Set $\varphi(t,\omega)x:=X_t(x)$, and then for an invariant Markov measure $\mu$
of the skew product flow $\Theta_t$, the stationary measure $\bar{\mu}$ for $p_t$ is given by
$\mE(\mu_{\cdot}|\mathscr{F}_{>0}).$
\el
\begin{proof}
By Definition \ref{Marm}, $\mE(\mu_{\cdot}|\mathscr{F}_{<\infty})=\mE(\mu_{\cdot}|\mathscr{F}_{=0})$.
Set $\bar{\mu}:=\mE(\mu_{\cdot}|\mathscr{F}_{>0})$, and then it follows from independence of $\mathscr{F}_{>0}$ and $\mathscr{F}_{=0}$ that
\ce
\bar{\mu}&=&\mE(\mu_{\cdot}|\mathscr{F}_{<\infty}\wedge\mathscr{F}_{>0})=\mE[\mE(\mu_{\cdot}|\mathscr{F}_{<\infty})|\mathscr{F}_{>0}]\\
&=&\mE\left[\mE(\mu_{\cdot}|\mathscr{F}_{=0})|\mathscr{F}_{>0}\right]=\mE\left[\mE(\mu_{\cdot}|\mathscr{F}_{=0})\right]
=\mE(\mu_{\cdot}).
\de
Therefore $\bar{\mu}$ is not random.

On one hand, by (\ref{inveq}), it holds that for $t>0$
\ce
\mE[\varphi(t,\cdot)\mu_{\cdot}|\theta_t^{-1}\mathscr{F}_{>0}]
&=&\mE[\mE(\varphi(t,\cdot)\mu_{\cdot}|\theta_t^{-1}\mathscr{F})|\theta_t^{-1}\mathscr{F}_{>0}]
=\mE[\mu_{\theta_t\cdot}|\theta_t^{-1}\mathscr{F}_{>0}]\\
&=&\mE(\mu_{\cdot}|\mathscr{F}_{>0})(\theta_t\omega)
=\bar{\mu}.
\de
On the other hand,
\ce
\mE[\varphi(t,\cdot)\mu_{\cdot}|\theta_t^{-1}\mathscr{F}_{>0}]
=\mE[\mE(\varphi(t,\cdot)\mu_{\cdot}|\mathscr{F}_{>0})|\theta_t^{-1}\mathscr{F}_{>0}]
=\mE[\varphi(t,\cdot)\bar{\mu}|\theta_t^{-1}\mathscr{F}_{>0}]
=\mE[\varphi(t,\cdot)\bar{\mu}].
\de
Thus for $B\in\sB(\mR^d)$,
\ce
\bar{\mu}(B)&=&\mE[\varphi(t,\cdot)\bar{\mu}(B)]=\mE[\bar{\mu}(\varphi(t,\cdot)^{-1}B)]
=\int_{\Omega}\mP(\dif \omega)\int_{\mR^d}1_B(\varphi(t,\omega)x)\bar{\mu}(\dif x)\\
&=&\int_{\mR^d}\bar{\mu}(\dif x)\int_{\Omega}1_B(\varphi(t,\omega)x)\mP(\dif \omega)
=\int_{\mR^d}p_t(x,B)\bar{\mu}(\dif x).
\de
By Definition \ref{stam} $\mE(\mu_{\cdot}|\mathscr{F}_{>0})$ is a stationary
measure for $p_t$.
\end{proof}

Next, under these conditions {\bf (H$_b$)} {\bf (H$_\sigma$)} {\bf (H$_f$)} {\bf (H$_g$)},
it follows from \cite[Theorem 1.3]{hqxz} that for almost all $\omega\in\Omega$, $x\mapsto X_t(x,\omega)\in\cH(\mR^d)$
for all $t\geq 0$. Set $\varphi(t,\omega)x:=X_t(x)$, and then $\varphi(t,\omega)$ is a
measurable c\`adl\`ag random dynamical system on the Polish space $(\mR^d,\sB(\mR^d))$.
When $(\mR^d,\sB(\mR^d))$ is compactificated as $(\hat{\mR}^{d},\sB(\hat{\mR}^{d}))$,
\cite[Lemma 5.1]{hc} admits us to obtain that an invariant Markov measure $\mu$
of the skew product flow $\Theta_t$ associated with $\varphi(t,\omega)$ exists. Set
$\bar{\mu}:=\mE(\mu_{\cdot}|\mathscr{F}_{>0})$, and then by Lemma \ref{invmeqstam} it holds
that $\bar{\mu}$ is a stationary measure for Eq.(\ref{Eqj1}). Thus, we get the following
theorem.

\bt\label{stae}
Suppose that these conditions {\bf (H$_b$)} {\bf (H$_\sigma$)} {\bf (H$_f$)} {\bf (H$_g$)} hold.
Then on the compact space $(\hat{\mR}^{d},\sB(\hat{\mR}^{d}))$, stationary measures for Eq.(\ref{Eqj1})
exist.
\et

\medskip

Next, we weaken {\bf (H$_f$)} {\bf (H$_g$)} and show by the similar Lyapunov approach to in \cite{bhan}
that stationary measures for Eq.(\ref{Eqj1}) still exist.

Assume that
\begin{enumerate}[{\bf (H$'_f$)}]
\item For some $q>(2d)\vee 4$ and any $p\in[2,q]$, there exists a constant $C_p>0$ such that
for $x,y\in\mR^{d}$
$$
\int_{\mU_{0}}|f(x,u)-f(y,u)|^{p}\,\nu(\dif u)\leq
C_p|x-y|^{p}\cdot\log(|x-y|^{-1}+\me),
$$
and
$$
\int_{\mU_{0}}|f(x,u)|^{p}\,\nu(\dif u)\leq C_p(1+|x|)^{p}.
$$
\end{enumerate}
\begin{enumerate}[{\bf (H$'_g$)}]
\item For each $u\in\mU-\mU_0$, $x\mapsto g(x,u)\in\cC(\mR^d)$.
\end{enumerate}

Under {\bf (H$_b$)}, {\bf (H$_\sigma$)}, {\bf (H$'_f$)} and {\bf (H$'_g$)}, it
is well known that there exists a unique strong solution to Eq.(\ref{Eqj1})(\cite{hqxz}).
This solution will be still denoted by $X_t(x)$. Moreover, the infinitesimal generator $\cL$ of $X_t(x)$
is given by
\ce
(\cL h)(y)&=&\<\partial_yh(y),b(y)\>
+\frac{1}{2}\frac{\partial^2 h(y)}{\partial y_i\partial y_j}\sigma_{ik}(y)\sigma_{kj}(y)\\
&&+\int_{\mU_0}\Big(h\big(y+f(y,u)\big)-h(y)-\<\partial_y h(y),f(y,u)\>\Big)\nu(\dif u)\\
&&+\int_{\mU-\mU_0}\Big(h\big(y+g(y,u)\big)-h(y)\Big)\nu(\dif u)
\de
for $h\in\cC^2_c(\mR^d)$ (\cite{bhan}). We introduce two key concepts.

\bd\label{ulbo}
$X_t(x)$ is said to be $2$-ultimately bounded if there exists a
positive finite constant $M$ such that
\ce
\limsup\limits_{t\rightarrow\infty}\mE|X_t(x)|^2\leq M.
\de
\ed

\bd\label{exulbo}
$X_t(x)$ is said to be exponentially $2$-ultimately bounded if
there exist positive finite constants $K, \beta, M$ such that
\ce
\mE|X_t(x)|^2\leq K\me^{-\beta t}|x|^2+M.
\de
\ed

Obviously, if $X_t(x)$ is exponentially $2$-ultimately bounded, then $X_t(x)$ is $2$-ultimately bounded.
Next, we show that $X_t(x)$ is exponentially $2$-ultimately bounded under some conditions.

\bp\label{laexul}
If there exists a function $V\in\cC^2(\mR^d)$ such that

(i) $K_1|x|^2-M_1\leq V(x)\leq K_2|x|^2+M_2$,

(ii) $\cL V(x)\leq -K_3V(x)+M_3$,\\
where $K_1>0, K_2>0, K_3>0, M_1\geq0, M_2\geq0, M_3\geq0$ are constants, then $X_t(x)$ is exponentially $2$-ultimately bounded.
\ep
\begin{proof}
Applying the It\^o formula to $e^{K_3t}V(X_t)$, we obtain that
\ce
e^{K_3t}V(X_t)-V(x)&=&\int_0^tK_3e^{K_3s}V(X_s)\dif s+\int_0^te^{K_3s}\<\partial_y V(X_s),b(X_s)\>\dif s\\
&&+\int_0^te^{K_3s}\<\partial_y V(X_s),\sigma(X_s)\dif W_s\>\\
&&+\int_0^t\int_{\mU_0}e^{K_3s}\(V\big(X_s+f(X_s,u)\big)-V(X_s)\)\tilde{N}_{\kappa}(\dif s, \dif u)\\
&&+\int_0^t\int_{\mU-\mU_0}e^{K_3s}\(V\big(X_s+g(X_s,u)\big)-V(X_s)\)N_{\kappa}(\dif s, \dif u)\\
&&+\frac{1}{2}\int_0^te^{K_3s}\frac{\partial^2V(X_s)}{\partial y_i\partial y_j}\sigma_{ik}(X_s)\sigma_{kj}(X_s)\dif s\\
&&+\int_0^t\int_{\mU_0}e^{K_3s}\Big(V\big(X_s+f(X_s,u)\big)-V(X_s)\\
&&-\<\partial_yV(X_s),f(X_s,u)\>\Big)
\nu(\dif u)\dif s.
\de
Taking expectation on two hands sides, we further get by (ii) that
\ce
e^{K_3t}\mE V(X_t)-V(x)&=&\int_0^tK_3e^{K_3s}\mE V(X_s)\dif s+\int_0^te^{K_3s}\mE\<\partial_y V(X_s),b(X_s)\>\dif s\\
&&+\int_0^t\int_{\mU-\mU_0}e^{K_3s}\mE\(V\big(X_s+g(X_s,u)\big)-V(X_s)\)\nu(\dif u)\dif s\\
&&+\frac{1}{2}\int_0^te^{K_3s}\mE\left[\frac{\partial^2V(X_s)}{\partial y_i\partial y_j}\sigma_{ik}(X_s)\sigma_{kj}(X_s)\right]\dif s\\
&&+\int_0^t\int_{\mU_0}e^{K_3s}\mE\Big(V\big(X_s+f(X_s,u)\big)-V(X_s)\\
&&-\<\partial_yV(X_s),f(X_s,u)\>\Big)
\nu(\dif u)\dif s\\
&=&\int_0^tK_3e^{K_3s}\mE V(X_s)\dif s+\int_0^te^{K_3s}\mE\cL V(X_s)\dif s\\
&\leq&\int_0^tK_3e^{K_3s}\mE V(X_s)\dif s+\int_0^te^{K_3s}\mE\left(-K_3V(X_s)+M_3\right)\dif s\\
&=&M_3\frac{e^{K_3t}-1}{K_3}.
\de
On one hand, it follows from (i) that
\be
\mE V(X_t)\leq e^{-K_3t}V(x)+M_3\frac{1-e^{-K_3t}}{K_3}
\leq e^{-K_3t}\left(K_2|x|^2+M_2\right)+M_3\frac{1-e^{-K_3t}}{K_3}.
\label{right}
\ee
On the other hand, (i) also admits us to have that
\be
\mE V(X_t)\geq K_1\mE|X_t|^2-M_1.
\label{left}
\ee
Combining (\ref{left}) with (\ref{right}), we obtain that
\ce
\mE|X_t|^2&\leq&\frac{K_2}{K_1}e^{-K_3t}|x|^2+\left(\frac{M_1}{K_1}+\frac{M_2e^{-K_3t}}{K_1}+M_3\frac{1-e^{-K_3t}}{K_3K_1}\right)\\
&\leq&\frac{K_2}{K_1}e^{-K_3t}|x|^2+\frac{K_3(M_1+M_2)+M_3}{K_3K_1}.
\de
The proof is completed.
\end{proof}

The function $V(x)$ in the above Proposition is called as a Lyapunov function of $X_t(x)$.
In the following, we use exponentially $2$-ultimately bounded property of $X_t(x)$ to prove existence
of stationary measures for Eq.(\ref{Eqj1}).

\bt\label{exi}
Suppose that there exists a Lyapunov function $V(x)$ of $X_t(x)$. Then Eq.(\ref{Eqj1}) has a stationary measure.
\et
\begin{proof}
By the above Proposition, we know that $X_t(x)$ is exponentially $2$-ultimately bounded and then
$2$-ultimately bounded. Thus by Definition \ref{ulbo}, there exists a $t_0>0$ such that for $t>t_0$,
$\mE|X_t(x)|^2<M$ and
\ce
\frac{1}{t}\int_0^t\mE|X_s|^2\dif s&=&\frac{1}{t}\left[\int_0^{t_0}\mE|X_s|^2\dif s+\int_{t_0}^t\mE|X_s|^2\dif s\right]\\
&\leq&\frac{1}{t}\left[\int_0^{t_0}\mE|X_s|^2\dif s+M(t-t_0)\right]\\
&=&\frac{1}{t}\int_0^{t_0}\mE|X_s|^2\dif s+M\left(1-\frac{t_0}{t}\right).
\de

Set
\ce
\bar{\mu}_T(B):=\frac{1}{T}\int_0^Tp_t(x,B)\dif t,
\de
for any $T>0$ and $B\in\mathscr{B}(\mR^d)$. And we have by Chebyshev's inequality that for $T>t_0$,
\ce
\bar{\mu}_T(B^c(0,R))&=&\frac{1}{T}\int_0^Tp_t(x,B^c(0,R))\dif t=\frac{1}{T}\int_0^T\mP(|X_t(x)|>R)\dif t\\
&\leq&\frac{1}{TR^2}\int_0^T\mE|X_t|^2\dif t
\leq\frac{1}{R^2}\left(\frac{1}{T}\int_0^{t_0}\mE|X_s|^2\dif s+M\left(1-\frac{t_0}{T}\right)\right).
\de
Thus, for any $\varepsilon>0$, $\bar{\mu}_T(B(0,R))>1-\varepsilon$ for
$R$ being large enough. Hence, $\{\bar{\mu}_T,T>t_0\}$ is tight and its
limit $\bar{\mu}$ is a stationary measure of Eq.(\ref{Eqj1}).
\end{proof}

\section{Special case 1: $f(x,u)=g(x,u)=u$}\label{fpe}

In the section, requiring that $f(x,u)=g(x,u)=u$ and $\nu(\dif u)=\frac{C_{d,\alpha}}{|u|^{d+\alpha}}\dif u$,
the L\'evy measure of a symmetric $\alpha$-stable
process, we show that stationary measures
for Eq.(\ref{Eqj1}) could be represented as the solutions for a type of Fokker-Planck equations.

Consider the following equation
\ce
X_{t}(x)=x+\int^{t}_{0}b(X_{s}(x))\,\dif s+\int^{t}_{0}\sigma(X_{s}(x))\,\dif W_{s}
+L_t,
\de
where $L_t$ is a symmetric $\alpha$-stable process independent
of $W_t$. Based on the l\'evy-It\^o representation of $L_t$, the above
equation can be rewritten as follows:

\begin{enumerate}[(i)]
\item for $1\leq\alpha<2$,
\be
X_{t}(x)&=&x+\int^{t}_{0}b(X_{s}(x))\,\dif s+\int^{t}_{0}\sigma(X_{s}(x))\,\dif W_{s}\no\\
&&+\int_0^t\int_{|u|\leq\delta}u\tilde{N}_\kappa(\dif s, \dif u)+\int_0^t\int_{|u|>\delta}uN_\kappa(\dif s, \dif u),
\label{al1}
\ee
where $\kappa_t:=L_t-L_{t-}$ and $\delta$ is the same one as in Section \ref{exst}.
\end{enumerate}

\begin{enumerate}[(ii)]
\item for $0<\alpha<1$,
\be
X_{t}(x)&=&x+\int^{t}_{0}b(X_{s}(x))\,\dif s+\int^{t}_{0}\sigma(X_{s}(x))\,\dif W_{s}\no\\
&&+\int_0^t\int_{\mR^d\setminus\{0\}}uN_\kappa(\dif s, \dif u).
\label{al0}
\ee
\end{enumerate}

We study mainly Eq.(\ref{al1}) and Eq.(\ref{al0}) can be dealt with similarly.

Under {\bf (H$_b$)} and {\bf (H$_\sigma$)}, by \cite[Theorem 1.3]{hqxz},
for almost all $\omega\in\Omega$, $x\mapsto X_t(x,\omega)$ is a homeomorphism
mapping on $\mR^d$, where $X_t(x,\omega)$ is the solution of Eq.(\ref{al1}).
Define
\ce
(p_th)(x):=\mE\big[h\big(X_t(x)\big)\big]=\int_{\mR^d}h(y)p_t(x,\dif y),
\de
for $h\in\cC_0(\mR^d)$. Thus $p_th\in\cC_0(\mR^d)$ by dominated convergence
theorem. Let $\cM_r(\mR^d)$ be the set of all finite regular signed measures
on $\sB(\mR^d)$. And then it is adjoint of $\cC_0(\mR^d)$ (\cite{ky}).

Next, we give some useful lemmas.

\bl\label{stco}
$X_t$ is stochastically continuous.
\el
\begin{proof}
By continuity of the integral and the stochastic integral, and the definition of stochastical continuity
(\cite[Definition 1.5, P.3]{sa}), one could obtain the result.
\end{proof}

\bl\label{sccs}
The family of operators $\{p_t, t\geq0\}$ defined above is a strongly continuous
contraction semigroup on $\cC_0(\mR^d)$.
\el
\begin{proof}
By \cite[Theorem 6.4.6, P.388]{da}, $\{p_t, t\geq0\}$ forms a semigroup.

Next, we prove strong continuity. For $h\in\cC_0(\mR^d)$, $h$ is uniformly
continuous on $\mR^d$. And for $\forall\varepsilon>0$, there exists an $\eta>0$
such that $|h(x)-h(y)|<\varepsilon$ for $x,y\in\mR^d$, $|x-y|<\eta$. For any
$\lambda\in\cM_r(\mR^d)$
\ce
&&\left|\int_{\mR^d}(p_th)(x)\lambda(\dif x)-\int_{\mR^d}h(x)\lambda(\dif x)\right|\\
&=&\left|\int_{\mR^d}\int_{\mR^d}h(y)p_t(x,\dif y)\lambda(\dif x)
-\int_{\mR^d}\int_{\mR^d}h(x)p_t(x,\dif y)\lambda(\dif x)\right|\\
&\leq&\int_{\mR^d}\int_{\mR^d}|h(y)-h(x)|p_t(x,\dif y)|\lambda|(\dif x)\\
&\leq&\int_{\mR^d}\int_{|x-y|<\eta}|h(y)-h(x)|p_t(x,\dif y)|\lambda|(\dif x)\\
&&+\int_{\mR^d}\int_{|x-y|\geq\eta}|h(y)-h(x)|p_t(x,\dif y)|\lambda|(\dif x)\\
&\leq&|\lambda|(\mR^d)\varepsilon+2\|h\|\int_{\mR^d}\mP\{|X_t(x)-x|\geq\eta\}|\lambda|(\dif x),
\de
where $|\lambda|$ stands for the variation measure of the signed measure
$\lambda$. For $\int_{\mR^d}\mP\{|X_t(x)-x|\geq\eta\}|\lambda|(\dif x)$,
by Lemma \ref{stco} and dominated convergence theorem,
when $t$ is small enough,
$$
\int_{\mR^d}\mP\{|X_t(x)-x|\geq\eta\}|\lambda|(\dif x)<\varepsilon.
$$
So,
\ce
\lim\limits_{t\downarrow0}\int_{\mR^d}(p_th)(x)\lambda(\dif x)=\int_{\mR^d}h(x)\lambda(\dif x).
\de
That is to say, $p_th$ converges weakly to $h$. By \cite[Theorem, p.233]{ky},
$p_th$ converges strongly to $h$.

Finally, by \cite[Definition 2.1, p.4]{ap}, $\{p_t, t\geq0\}$ is a strongly
continuous contraction semigroup on $\cC_0(\mR^d)$.
\end{proof}

Let $\cL$ be the infinitesimal generator of $\{p_t, t\geq0\}$ and then for $h\in\cC^2_c(\mR^d)$ (\cite{da})
\ce
(\cL h)(y)=\<\partial_y h(y),b(y)\>
+\frac{1}{2}\frac{\partial^2 h(y)}{\partial y_i\partial y_j}\sigma_{ik}(y)\sigma_{kj}(y)
+(\cL_{\alpha} h)(y).
\de
By Theorem \ref{sead}, we know that the adjoint $\cL^*$ of $\cL$ is given by
$$
(\cL^*\psi)(y):=-\frac{\partial}{\partial y_j}\big(b_j(y)\psi(y)\big)
+\frac{1}{2}\frac{\partial^2}{\partial y_i\partial y_j}\big(\sigma_{ik}(y)\sigma_{kj}(y)\psi(y)\big)
+(\cL_\alpha \psi)(y)
$$
for $\psi\in\cC^\infty_c(\mR^d)$. The main result in this section is the following theorem.

\bt\label{fkst}
If there exists a $\rho\in\cD(\cL^*)$ satisfying the following equation
\ce
\cL^*\rho=0
\de
with these conditions $\rho(y)\geq 0, \forall y\in\mR^{d}$ and $\int_{\mR^{d}}\rho(y)\dif y=1$.
Then $\bar{\mu}(\dif y):=\rho(y)\dif y$ is a stationary measure for $p_t$.
\et
\begin{proof}
By Lemma \ref{sccs}, $\{p_t, t\geq0\}$ is a strongly continuous
semigroup on $\cC_0(\mR^d)$ with the infinitesimal generator $\cL$. Let $p_t^*$
be adjoint of $p_t$ and $\overline{\cD(\cL^*)}$ stand for the closure of $\cD(\cL^*)$ in $\cM_r(\mR^d)$.
\cite[Theorem 10.4, p.41]{ap} admits us to get
that the restriction $p_t^+$ of $p_t^*$ to $\overline{\cD(\cL^*)}$ is a strongly
continuous semigroup on $\cM_r(\mR^d)$. Moreover, the infinitesimal
generator $\cL^+$ of $p_t^+$ is the part of $\cL^*$ in $\overline{\cD(\cL^*)}$,
i.e. $\cD(\cL^+)=\{h\in\cD(\cL^*), \cL^*h\in\overline{\cD(\cL^*)}\}$
and $\cL^+h=\cL^*h$ for $h\in\cD(\cL^+)$.

Since $\rho\in\cD(\cL^*)$ and $\cL^*\rho=0$, $\rho\in\cD(\cL^+)$ and $\cL^+\rho=0$.
By \cite[Theorem 2.4, p.4]{ap}, we obtain that
$$
p_t^+\rho=\rho, \quad t\geq0.
$$
So, for $\varphi\in\cC_c(\mR^d)\subset\cC_0(\mR^d)$,
\ce
\int_{\mR^{d}}\varphi(x)\rho(x)\dif x=\int_{\mR^{d}}\varphi(x)(p_t^+\rho)(x)\dif x=\int_{\mR^{d}}\varphi(x)(p_t^*\rho)(x)\dif x
=\int_{\mR^{d}}(p_t \varphi)(x)\rho(x)\dif x.
\de
Besides, under these conditions $\rho(y)\geq 0, \forall y\in\mR^{d}$ and $\int_{\mR^{d}}\rho(y)\dif y=1$,
$\rho(y)$ is a density function and $\bar{\mu}(\dif y):=\rho(y)\dif y$ is a probability measure.
Thus, by density of $\cC_c(\mR^d)$ in $L^2(\mR^d,\cB(\mR^d),\bar{\mu})$, it holds that for $B\in\cB(\mR^d)$
\ce
\int_{\mR^{d}}p_t(x,B)\bar{\mu}(\dif x)=\int_{\mR^{d}}(p_t I_B)(x)\bar{\mu}(\dif x)=\int_{\mR^{d}}I_B(x)\bar{\mu}(\dif x)=\bar{\mu}(B).
\de
By Definition \ref{stam}, $\bar{\mu}(\dif y)$ is a stationary measure
for $p_t$.
\end{proof}

\br\label{fksm}
If $b(x)=-x$ and $\sigma(x)=0$, the above theorem is \cite[Proposition 3.2(ii)]{arw}.
Moreover,
$$
\hat{\rho}(u)=\exp\{-\frac{1}{\alpha}C|u|^\alpha\}, \quad u\in\mR^{d},
$$
where $C$ is the same constant as one in Definition \ref{rid}.
\er

\section{Special case 2: $g(x,u)=0$}\label{des}

In the section, we require that $g(x,u)=0$ in Eq.(\ref{Eqj1}). And then Eq.(\ref{Eqj1})
is changed as the following equation
\be
X_{t}(x)&=&x+\int^{t}_{0}b(X_{s}(x))\,\dif s+\int^{t}_{0}\sigma(X_{s}(x))\,\dif W_{s}\no\\
&&+\int^{t+}_{0}\int_{\mU_{0}}f(X_{s-}(x),u)\,\tilde{N}_{k}(\dif s,\dif u), \quad t\geq0.
\label{Eqj2}
\ee
For Eq.(\ref{Eqj2}), if $b, \sigma, f$ satisfy some regular conditions, we could
know what stationary measures for Eq.(\ref{Eqj2}) are certainly.

{\bf Assume:}

\begin{enumerate}[({\bf H}$^1_{b,\sigma,f}$)]
\item $b$ and $\sigma$ are $(4d+6)$-times differentiable with bounded
derivatives of all order between $1$ and $4d+6$. Besides, $f(\cdot,u)$
is $(4d+6)$-times differentiable, and
\ce
f(0,\cdot)&\in&\bigcap_{2\leq q<\infty}L^q(\mU_0,n)\\
\sup\limits_{x}|\partial^r_xf(x,\cdot)|&\in&\bigcap_{2\leq q<\infty}L^q(\mU_0,n),
\quad 1\leq r\leq4d+6,
\de
where the space $(\mU_0,n)$ is equipped with a norm and
$\partial^r_xf(x,\cdot)$ stands for $r$ order partial derivative of $f(x,\cdot)$
with respect to $x$.
\end{enumerate}

\begin{enumerate}[({\bf H}$^2_{b,\sigma,f}$)]
\item There exist three constants $\varepsilon>0$, $\eta\geq0$ and $C>0$ such that
for all $x,y\in\mR^{d}$
\ce
\<y,\sigma(x)\sigma^T(x)y\>\geq|y|^2\frac{\varepsilon}{1+|x|^\eta},
\qquad |\det\{I+r\partial_xf(x,u)\}|\geq C ~\mbox{for all}~r\in[0,1].
\de
\end{enumerate}

Under ({\bf H}$^1_{b,\sigma,f}$) and ({\bf H}$^2_{b,\sigma,f}$), by
\cite[Theorem 2-29, p.15]{bgj}, Eq.(\ref{Eqj2}) has a unique solution
denoted by $X_t(x)$. Moreover, the transition probability $p_t(x,\dif y)$
has a density $\rho_t(x,y)$ and $(t,x,y)\mapsto\rho_t(x,y)$ is continuous.
Thus, the distribution of $X_t$ has the density $\rho_t(x,y)$.

\bt\label{list}
Suppose that $\lim\limits_{t\rightarrow\infty}\rho_t(x,y)=\rho(y)$, where
$\rho(y)$ satisfies
\ce
\rho(y)\geq 0, \quad \forall y\in\mR^{d} \quad \mbox{and} \quad \int_{\mR^{d}}\rho(y)\dif y=1.
\de
Then $\bar{\mu}(\dif y):=\rho(y)\dif y$ is a stationary measure for $p_t$.
\et
\begin{proof}
For $t>0$ and $B\in\cB(\mR^d)$,
\ce
\int_{\mR^{d}}p_t(x,B)\bar{\mu}(\dif x)&=&\int_{\mR^{d}}p_t(x,B)\rho(x)\dif x
=\int_{\mR^{d}}p_t(x,B)\lim\limits_{s\rightarrow\infty}\rho_s(x,x)\dif x\\
&=&\lim\limits_{s\rightarrow\infty}\int_{\mR^{d}}p_t(x,B)\rho_s(x,x)\dif x
=\lim\limits_{s\rightarrow\infty}\int_{\mR^{d}}\int_B\rho_t(x,y)\rho_s(x,x)\dif y\dif x\\
&=&\lim\limits_{s\rightarrow\infty}\int_B\int_{\mR^{d}}\rho_t(x,y)\rho_s(x,x)\dif x\dif y
=\lim\limits_{s\rightarrow\infty}\int_B\rho_{s+t}(x,y)\dif y
=\bar{\mu}(B).
\de
By Definition \ref{stam}, $\bar{\mu}(\dif y)$ is a stationary measure
for $p_t$.
\end{proof}

\br\label{lim}
By the above theorem, we see that if a limiting distribution exists, it must be a
stationary measure. This theorem also has a corresponding version in the theory of
Markov chains (\cite[p.237]{mt}).
\er

\end{document}